\documentclass[twoside,11pt]{amsart}
\usepackage{amsmath,latexsym,amssymb,times,mathptm, enumerate}

\def\proof{\noindent{\bf{Proof.} }}
\def\sqr#1#2{{\vcenter{\hrule height.#2pt
        \hbox{\vrule width.#2pt height#1pt \kern#1pt
                \vrule width.#2pt}
        \hrule height.#2pt}}}
\def\square{\mathchoice\sqr64\sqr64\sqr{4}3\sqr{3}3}
\def\QED{\hfill$\square$}

\def\tratto{\mbox{\rule{2mm}{.2mm}$\;\!$}}

\newtheorem{theorem}{Theorem}[section]

\newtheorem{corollary}[theorem]{Corollary}
\newtheorem{lemma}[theorem]{Lemma}
\newtheorem{proposition}[theorem]{Proposition}

 \theoremstyle{definition}

\newtheorem{remark}[theorem]{Remark}

\newtheorem{example}[theorem]{Example}

\def\p{{\mathfrak p}}
\def\a{{\mathfrak a}}
\def\A{{\mathfrak A}}
\newcommand{\s}{\bigskip}
\newcommand{\ms}{\medskip}

\newcommand{\f}[1]{\ensuremath{\mathfrak{#1}}}
\newcommand{\ol}[1]{\ensuremath{\overline{#1}}}
\newcommand{\mc}[1]{\ensuremath{\mathcal{#1}}}
\def\m{{\mathfrak m}}

\numberwithin{equation}{section}

\def\noi{\noindent }

\newcommand{\core}[1]{\ensuremath{{\rm{core}}(#1)}}
\begin{document}
\baselineskip=16pt

\title{The Core of Ideals in Arbitrary Characteristic}
\author[L. Fouli]{Louiza Fouli}
\address{Department of Mathematics, University of Texas, Austin, Texas 78712, USA}
\email{lfouli@math.utexas.edu}
\author[C. Polini]{Claudia Polini}
\address{Department of Mathematics, University of Notre Dame
, Notre Dame, Indiana 46556, USA} \email{cpolini@nd.edu}
\author[B. Ulrich]{Bernd Ulrich}
\address{Department of Mathematics, Purdue University
, West Lafayette, Indiana 47907, USA} \email{ulrich@math.purdue.edu}


\thanks{The second and third author were supported in part by the
NSF. The second author was also supported in part by the NSA}

\vspace{-0.1in}

\vspace{-0.2in}

\maketitle

\section{Introduction}

In this article we provide explicit formulas for the core of an
ideal. Recall that for  an ideal $I$ in a Noetherian ring $R$, the
{\it core} of $I\/$, $\core{I}$, is the intersection of all
reductions of $I\/$. For a subideal $J \subset I$ we say $J\/$ is a
{\it reduction\/} of $I\/$, or  $I\/$ is {\it integral } over $J\/$,
if $I^{r+1}=JI^{r}$ for some $r \geq 0\/$; the smallest such $r$ is
called the {\it reduction number} of $I\/$ with respect to $J\/$ and
is denoted by $r_J(I)$. If $(R, \m)$ is local with infinite residue
field $k$, every ideal has a {\it minimal reduction}, which is a
reduction minimal with respect to inclusion. Minimal reductions of a
given ideal $I\/$ are far from being unique, but they all share the
same minimal number of generators, called the {\it analytic spread}
of $I\/$ and written $\ell(I)$. Minimal reductions arise from
Noether normalizations of the {\it special fiber ring \,}
$\mc{F}(I)= {\rm gr}_I(R) \otimes k$ of $I\/$, and therefore
$\ell(I)={\rm{dim}} \, \mc{F}(I)$. From this one readily sees that ${\rm{ht}} \,  I
\leq \ell(I) \leq {\rm{dim}} \, R$; these inequalities are equalities for any
$\m$-primary ideal, and if the first inequality is an equality then
$I\/$ is called {\it equimultiple}. Obviously, the core can be
obtained as an intersection of minimal reductions of a given ideal.

Through the study of the core one hopes to better understand
properties shared by all reductions. The notion was introduced by
Rees and Sally with the purpose of generalizing the
Brian\c{c}on-Skoda Theorem \cite{RS}. Being an a priori infinite
intersection of reductions the core is difficult to compute, and
there have been considerable efforts to find explicit formulas, see
\cite{HS, CPU01, CPU02, HySm1, PU, HT, HySm2}. We quote the
following result from \cite{PU}:
\begin{theorem}\label{PUmain}
Let $R$ be a local Gorenstein ring with infinite residue field $k$,
let $I\/$ be an $R$-ideal with $g={\rm{ht}} \, I > 0$ and $\ell=
\ell(I)$, and let $J\/$ be a minimal reduction of $I\/$ with
$r=r_J(I)$. Assume $I\/$ satisfies $G_{\ell}$, ${\rm depth} \ R/I^j
\geq {\rm dim} \ R/I -j+1$ for $1 \leq j \leq \ell-g$, and either
${\rm char} \ k =0 $ or ${\rm char} \ k  > r -\ell + g $. Then
\[
{\rm core}(I)= J^{\/n+1} : I^{\/n}
\]
for every $n \geq {\rm max} \{ r - \ell + g, 0 \}$.
\end{theorem}

The property $G_{\ell}$ in the above theorem is a rather weak
requirement on the local number of generators of $I\/$: It means
that the minimal number of generators $\mu(I_{\mathfrak p})$ is at
most ${\rm{dim}} \, R_{\p}$ for every prime ideal ${\mathfrak p}$ containing
$I\/$ with ${\rm{dim}} \, R_{\p}\leq \ell -1 $. Both hypotheses, the
$G_{\ell}$ condition and the depth assumption on the powers, are
automatically satisfied if $I\/$ is equimultiple. They also hold for
one-dimensional generic complete intersection ideals, or more
generally, for Cohen-Macaulay generic complete intersections with
$\ell=g+1$. In the presence of the $G_{\ell}$ property, the depth
inequalities for the powers hold if $I\/$ is perfect with $g=2$,
$I\/$ is perfect Gorenstein with $g=3$, or more generally, if $I\/$
is in the linkage class of a complete intersection  \cite[1.11]{H1}.

Theorem~\ref{PUmain} is not true in general without the assumption
on the characteristic, as was shown in \cite[4.9]{PU}. Hence in the
present paper we study the case of arbitrary characteristic.
Explicit formulas for the core valid in any characteristic and for
any reduction number are known for equimultiple ideals of height one
\cite[3.4$($a$)$]{PU} and for powers of the homogeneous maximal
ideal of standard graded reduced Cohen-Macaulay rings over an
infinite perfect field \cite[4.1]{HySm2}. In this paper we clarify
the latter result and generalize it to ideals generated by forms of
the same degree that are not necessarily zero-dimensional or even
equimultiple:

\begin{theorem}\label{formsofthesamedegree}
Let $k$ be an infinite field, $R'$ a positively graded geometrically
reduced Cohen-Macaulay $k$-algebra, and $R$ the localization of $R'$
at the homogeneous maximal ideal. Let $I\/$ be an $R$-ideal
generated by forms in $R'$ of the same degree with $g={\rm ht} \,
{I} >0$ and $\ell=\ell(I)$, and let $J\/$ be a minimal reduction of
$I\/$ with $r=r_{J}(I)$. If $\ell > g$ further assume that $R'$ is
Gorenstein, $I\/$ satisfies $G_{\ell}$ and ${\rm{depth}} \ R/I^{j}
\geq {\rm{dim}} \ R/I-j+1$ for $1\leq j \leq \ell-g$. Then
\[
{\rm{core}}(I)= J^{n+1}:I^{n}
\]
for every $n \geq {\rm{max}} \{\;r-\ell+g,0\; \}$.

\end{theorem}

Recall that the $k$-algebra $R'$ is said to be {\it geometrically
reduced} if after tensoring with the algebraic closure ${\overline k}$
of $k$ the ring $R' \otimes_k {\overline k}$ is reduced.

The above theorem is a special case of a considerably more general
result in which the assumption on the grading is replaced by the
condition that the residue field is perfect and the special fiber
ring $\mc{F}(I)$ is reduced, or yet more generally, has embedding
dimension at most one locally at every minimal prime of maximal
dimension (Theorem~\ref{emdimformula}). We identify further
instances where the assumption on the special fiber ring is
satisfied: More generally than in Theorem~\ref{formsofthesamedegree}
it suffices to require that $I=(K,f)$ where $K$ is generated by
forms of the same degree and either $f$ is integral over $K$
(Theorem~\ref{cor integral}) or else $\ell(K) \leq \ell(I)$ and
$\mc{F}(K)$ satisfies Serre's condition $R_1$
(Theorem~\ref{regular}). We give a series of examples showing that
our hypotheses are sharp: Theorem~\ref{formsofthesamedegree} fails
to hold without the assumption of geometric reducedness, even when
$R'$ is a domain and $I =\m$ (Example~\ref{reduced}); this also
shows that an assumption needs to be added in \cite[4.1]{HySm2}.
Likewise, in Theorem~\ref{emdimformula} it does not suffice to
suppose that the generic embedding dimension of $\mc{F}(I)$ be at
most two (Example~\ref{emb2}), and in Theorems~\ref{cor integral}
and \ref{regular} we must require $f$ to be integral over $K$ or
$\mc{F}(K)$ to satisfy $R_1$ (Example~\ref{serre}).

Our approach, which is different from the one of \cite{HySm2}, can
be outlined as follows: Write $\ell=\ell(I)$, let $f_1, \ldots,
f_{\ell+1}$ be general elements in $I$, set $J=(f_1, \ldots,
f_{\ell})$, and let $^{\tratto}$ denote reduction modulo the
`geometric residual intersection' $(f_1, \ldots, f_{\ell -1}) : I$.
As $\overline{I}$ is an equimultiple ideal of height one, we can
apply the formula of \cite[3.4$($a$)$]{PU}, which says that
regardless of characteristic,
\[
{\rm core} (\ol{I}) =  \ol{J}^{n+1} : \displaystyle \sum_{y \in
\ol{I}} (\ol{J}, y)^n \hspace{1cm} \mbox{for} \ n\gg 0\,.
\]
The problem is that this formula does not `lift' from ${\overline
I}\/$ to $I\/$. On the other hand, according to one of our main
technical results, the equality $ {\rm core} ({\overline I}) =
\ol{J}^{n+1} : {\overline I}^n$ does lift (Lemma~\ref{abars}, see
also \cite[4.2]{PU}). Thus the task becomes to show that
\[
\displaystyle \sum_{y \in \ol{I}} (\ol{J}, y)^n= \ol{I}^n
\hspace{1cm} \mbox{for} \ n\gg 0\,.
\]
This follows from a general `decomposition formula' for powers that
may be interesting in its own right: In fact we prove that if $R$ is
a Noetherian local ring with infinite perfect residue field  and
$\mc{F}(I)$ has embedding dimension at most one locally at every
minimal prime of maximal dimension, then
\[
I^n = (f_1, \ldots, f_{\ell -1})I^{n-1} + (f_{\ell}, f_{\ell +1})^n
\mbox{ \ \ \ for \, } n\gg 0 \,
\]
(special case of Theorem~\ref{I^n decomposition}).

\s

\s

\section{A Decomposition Formula for Powers}

\s

In this section we show our decomposition formula for powers of
ideals. The proof is based on Theorem \ref{existence}, a generalization of the Primitive
Element Theorem. We begin by reviewing two lemmas:


\begin{lemma}\label{Lemma1}
Let $k$ be an infinite field, $A=k[X_1, \ldots, X_n]$ a polynomial
ring with quotient field $K$,
and $B$ an $A$-algebra essentially of finite type. Then \[ {\rm dim}
\ B \otimes_{A} A/(\{X_i - \lambda_i\}) \leq {\rm dim} \ B
\otimes_{A} K
\]
for $(\lambda_1, \ldots, \lambda_n) \in k^n$ general.
\end{lemma}
\proof By the Generic Flatness Lemma there exists an element
$0\not=f \in A$ so that $A_{f} \longrightarrow B_{f}$ is flat and
hence satisfies going down \cite[24.1]{M}. For every $(\lambda_1,
\ldots, \lambda_n) \in k^n \setminus V(f)$ one has
 \[ {\rm dim}
\ B \otimes_{A} A/(\{X_i - \lambda_i\}) \leq {\rm dim} \ B
\otimes_{A} K
\]
according to \cite[15.3]{M}. \QED

\s

\begin{lemma}\label{Lemma2}
Let $k$ be an infinite field, $C$ a finitely generated $k$-algebra,
and $I=(f_1, \ldots, f_n)$ a $C$-ideal.  Let $\a$ be a $C$-ideal
generated by $t$ general $k$-linear combinations of $f_1, \ldots,
f_n$. Then
\[ {\rm dim} \ C / (\a : I^{\infty})  \leq {\rm dim} \, C -t.
\]
In particular ${\rm dim} \ C / \a  \leq {\rm max} \{ {\rm dim} \, C
-t , {\rm dim} \ C/I\}$.
\end{lemma}

\proof Let $X_{ij}$ be variables over $k$, where $1 \leq i \leq t$
and $1 \leq j \leq n$, set $R = C[\{X_{ij}\}]$, and write $\A$ for
the $R$-ideal generated by the $t$ generic linear combinations
$\sum_{j=1}^{n} X_{ij}f_j$, where $1 \leq i \leq t$. We first show
that $\A : I^{\infty}$ has height at least $t$ in $R$, or
equivalently, that $IR \subset \sqrt{\A}$ locally in codimension at
most $t-1$. Thus let $Q$ be a prime ideal of $R$ that has height at
most $t-1$ and does not contain $I$. Replacing $C$ by $C_{Q \cap C}$
we may assume that $C$ is local and $I=C$, and after applying a
$C$-automorphism of $R$ we are in the situation where $f_1, \ldots,
f_n= 1,0, \ldots, 0$. But then $\A= (X_{11}, \ldots, X_{t1})$, which
cannot be contained in $Q$ as $Q$ has height at most $t-1$.

Next, consider the map $A = k[\{X_{ij}\}] \longrightarrow B=R/\A :
I^{\infty}$. Write $K$ for the quotient field of $A$, and $S=R
\otimes_A K =C \otimes_k K$. Notice that  ${\rm dim} \, S = {\rm dim}
\ C \otimes_{k} K ={\rm dim} \, C$ as $C$ is a finitely generated
$k$-algebra, and that ${\rm{ht}}(\A :
I^{\infty})S \geq {\rm ht} (\A : I^{\infty}) \geq t$.
Therefore
\[{\rm dim} \ B \otimes_A K = {\rm dim} \ S/(\A :I^{\infty})S
\leq {\rm dim} \, S - {\rm ht} (\A : I^{\infty})S \leq {\rm dim} \, C - t \, . \]
Finally, for a point
$(\lambda_{ij}) \in k^{tn}$ let $\a$ denote the $C$-ideal generated
by the $t$ elements $\sum_{j=1}^{n} \lambda_{ij}f_j$. Observe that $B
\otimes_{A} A/(\{X_{ij} - \lambda_{ij}\})$ maps onto  $C / (\a :
I^{\infty})$. Hence Lemma~\ref{Lemma1} shows that if
$(\lambda_{ij})$ is general then ${\rm dim} \ C / (\a : I^{\infty})
\leq {\rm dim} \ B \otimes_A K \leq {\rm dim} \, C -t$. \QED

\s

\begin{theorem}{\label{existence}}
Let $k$ be an infinite perfect field, $B=k[y_1, \ldots ,y_n]$ a
finitely generated $k$-algebra of dimension $d$, and $s$ a positive
integer. Let $A$ be a $k$-subalgebra generated by $d+s$ general
$k$-linear combinations of $y_1, \ldots, y_n$. Then $B$ is a finite
$A$-module, and ${\rm{dim}}_{A} \, B/A < d$ if and only if $B$ has
embedding dimension at most $s$ locally at every minimal prime of
dimension $d$.
\end{theorem}
\proof Clearly $B$ is a finite $A$-module by Lemma~\ref{Lemma2}.

First assume that ${\rm{dim}}_{A} \, B/A < d$. Let $\f{q} \in
{\rm{Spec}}(B)$ with ${\rm{dim}}\ B/\f{q}=d$. Let $\f{p}=\f{q} \cap
A$. Notice that ${\rm{dim}} \ A/\f{p}=d$. Since ${\rm{dim}}_{A} \,
B/A < d$ we have $A_{\f{p}}=B_{\f{p}}=B_{\f{q}}$. Write $A$ as an
epimorphic image of the polynomial ring $k[X_{1}, \ldots ,X_{d+s}]$
and let $\f{P}$ be the preimage of $\f{p}$ in $k[X_1, \ldots
,X_{d+s}]$. One has ${\rm{dim}} \ k[X_1, \ldots
,X_{d+s}]_{\f{P}}=d+s-{\rm{dim}} \ {k[X_1, \ldots
,X_{d+s}]}/\f{P}=d+s- \dim \ A/\f{p}=s$. Hence
$B_{\f{q}}=A_{\f{p}}$ has embedding dimension at most $s$.

We now assume that $B$ has embedding dimension at most $s$ locally
at every minimal prime of dimension $d$. Let $x_{1}, \ldots,
x_{d+s}$ be general $k$-linear combinations of $y_{1}, \ldots ,
y_{n}\/$. Consider the exact sequence
\[
0 \longrightarrow {\mathbb D} \longrightarrow C=B \otimes_{k}B
\overset{{\rm mult}}{\longrightarrow} B \longrightarrow 0 \ .
\]
\noi Notice that $\Omega_{k}(B)=\mathbb{D}/{\mathbb{D}}^{2}$ is the
module of differentials of $B$ over $k$. The $C$-ideal ${\mathbb D}$
is generated by $c_{i}=y_{i} \otimes 1 -1 \otimes y_{i}$ for $1 \leq
i \leq n$. Thus setting $a_{i}=x_{i} \otimes 1-1 \otimes x_{i}$ we
have that $a_{1}, \ldots , a_{d+s}$ are general $k$-linear
combinations of the generators $c_{1}, \ldots , c_{n}$ of
${\mathbb{D}}$. Write $\f{a}$ for the $C$-ideal generated by $a_{1},
\ldots , a_{d+s}$. According to Lemma~\ref{Lemma2}  \[ {\rm dim} \
C/ (\a : {\mathbb{D}}^{\infty})  \leq {\rm dim} \, C -d -s \leq  d-1.
\]
Hence for every $Q \in {\rm{Spec}} (C)$ with
${\rm{dim}} \ C/Q =d$ one has
$Q \in V(\f{a})$ if and only if
$Q \in V({\mathbb{D}})$. Observe that there are only finitely many
such primes as they are all minimal over ${\mathbb{D}}$.

Let $Q$ be one of these primes and write
$\f{q}=Q/{\mathbb{D}}$. Now ${\rm{dim}} \ B/\f{q} = d$ and hence
${\rm{edim}} \, B_{\f{q}} \leq s$. Consider the exact sequence
\[
\f{q}B_{\f{q}}/\f{q}^2B_{\f{q}} \longrightarrow
\Omega_{k}(B_{\f{q}}) \otimes_{B_{\f{q}}}k(\f{q}) \longrightarrow
\Omega_{k}(k(\f{q})) \longrightarrow 0 \ .
\]

\noi In this sequence we have $\mu(\f{q}B_{\f{q}}/\f{q}^2B_{\f{q}})
\leq s$, and $\mu(\Omega_{k}(k(\f{q})))= {\rm{trdeg}}_{k} \,
k(\f{q})= {\rm{dim}} \ B/\f{q} =d$ as $k$ is perfect. Thus
$\mu(\Omega_{k}(B_{\f{q}}) \otimes_{B_{\f{q}}}k(\f{q})) \leq d+s$.
Notice that
$\Omega_{k}(B_{\f{q}})={\mathbb{D}}_{Q}/{\mathbb{D}}_{Q}^{2}$, and
hence $\mu({\mathbb{D}}_{Q}) \leq d+s$ by Nakayama's Lemma. Therefore
${\mathbb{D}}_{Q}=\f{a}_{Q}$ by the general choice of $a_{1}, \ldots
, a_{d+s}$. In summary we obtain ${\mathbb{D}}_{Q}= \f{a}_{Q}$ for
every $Q \in V(\f{a})$ with ${\rm{dim}} \ C/Q=d$.

Write $A=k[x_{1}, \ldots ,x_{d+s}]$ and consider the exact sequence
\[
0 \longrightarrow {\mathbb{D}}'={\mathbb{D}}/ \f{a} \longrightarrow
C'=C/ \f{a}=B \otimes_{A} B \longrightarrow B \longrightarrow 0 \ .
\]
By the above ${\mathbb{ D}}'_{Q}=0$ for every $Q \in \rm{Spec}(C')$
with ${\rm{dim}} \ C'/Q=d$. The homomorphism $A \rightarrow C'=B
\otimes_{A} B$ makes $C'$ a finite $A$-module. Let $\f{p} \in
\rm{Spec}(A)$ with ${\rm{dim}} \ A/\f{p}=d$. Let $Q$ be any prime of
$C'$ lying over $\f{p}$. As $ {\rm{dim}} \ C'/Q= {\rm{dim}} \
A/\f{p}=d$ we obtain ${\mathbb{D}}'_{Q}=0$. Since this holds for any
such $Q$ we have ${\mathbb{D}}'_{\f{p}}=0$. Thus $B_{\f{p}}
\otimes_{A_{\f{p}}} B_{\f{p}} \cong B_{\f{p}}$. Computing numbers of
generators as $A_{\f{p}}$-modules we conclude $B_{\f{p}}=
A_{\f{p}}$, hence $(B/A)_{\f{p}}=0$.\QED

\s

Using the above result we are able to prove various versions of our
decomposition formula:

\begin{lemma}\label{powers}
Let $k$ be a field and $B$ a standard graded $k$-algebra of
dimension $d$ with homogeneous maximal ideal $\f{m}$. Let $A$ be a
$k$-subalgebra generated by $d+s$ linear forms $x_{1}, \ldots,
x_{d+s}$ in $B$. Then $\f{m}^{n}=(x_{1}, \ldots,
x_{d-1})\f{m}^{n-1}+(x_{d}, \ldots, x_{d+s})^{n}$ for $n \gg 0$ if
and only if $B/A$ is a finite module over $k[x_{1}, \ldots,
x_{d-1}]$ .
\end{lemma}
\proof Write $C=B/A$. Mapping variables $X_{i} \mapsto x_{i}$, we
obtain homogeneous maps
\[
k[X_{d}, \ldots, X_{d+s}] \twoheadrightarrow A/(x_{1}, \ldots,
x_{d-1})A \rightarrow B/ (x_{1}, \ldots, x_{d-1})B.
\]
Their composition is surjective in large degrees if and only if
$ C/(x_{1},\ldots,x_{d-1})C$ is a finite dimensional $k$-vector
space, which by the graded Nakayama Lemma means that $C$ is a finite
module over $k[x_{1}, \ldots, x_{d-1}]$ . \QED

\s

\begin{proposition}\label{maxidealpowers}
Let $k$ be an infinite perfect field, $B$ a standard graded
$k$-algebra of dimension $d$ with homogeneous maximal ideal $\f{m}$,
and $s$ a positive integer. Let $A$ be a $k$-subalgebra generated by
$d+s$ general linear forms $x_{1}, \ldots , x_{d+s}$ in $B$. Then
$\f{m}^{n}=(x_{1}, \ldots, x_{d-1})\f{m}^{n-1}+(x_{d}, \ldots,
x_{d+s})^{n}$ for $n \gg 0$ if and only if $B$ has embedding
dimension at most $s$ locally at every minimal prime of dimension
$d$.
\end{proposition}
\proof The assertion is an immediate consequence of
Theorem~\ref{existence} and Lemma~\ref{powers}.  \QED

\s

\begin{corollary} \label{finitecoker}
Let $R$ be a Noetherian local ring and $I\/$ an $R$-ideal of
analytic spread $\ell$. Let $f_{1}, \ldots , f_{\ell+s}$ be elements
in $I\/$, $\f{a}=(f_{1}, \ldots, f_{\ell-1})$, $K=(f_{1}, \ldots,
f_{\ell+s})$, and consider the natural map of special fiber rings
$\varphi :\mathcal{F}(K) \longrightarrow \mathcal{F}(I)$. Then
$I^{n}=(f_{1}, \ldots , f_{\ell-1})I^{n-1}+(f_{\ell}, \ldots,
f_{\ell+s})^{n}$ for $n \gg 0$ if and only if
${\rm{coker}}(\varphi)$ is a finite $\mathcal{F}(\f{a})$-module.
\end{corollary}
\proof We apply Lemma \ref{powers} with $B=\mathcal{F}(I)$ and
$A=\varphi(\mathcal{F}(K))$, and use Nakayama's Lemma.  \QED

\s

We are now ready to prove the main result of this section. Let $I$
be an ideal in a Noetherian local ring $R$ with infinite residue
field $k$.
Elements $f_1, \ldots, f_t$ in $I$ are
said to be {\it general} if the image of the tuple $(f_1, \ldots, f_t)$
is  a general point of the
affine space $(I \otimes k)^t$.  Recall that $t \geq \ell(I)$
general elements in $I$ generate a reduction,
and hence give  $I^n=(f_{1}, \ldots, f_t)I^{n-1}$ for $n \gg 0$.
The next result provides, under suitable assumptions, a different
type of decomposition formula for the powers of $I\,$:

\begin{theorem}\label{I^n decomposition}
Let $R$ be a Noetherian local ring with infinite perfect residue
field, $I\/$ an $R$-ideal of analytic spread $\ell$, and $s$ a
positive integer. Let $f_{1}, \ldots , f_{\ell+s}$ be general
elements in $I\/$, $\f{a}=(f_{1}, \ldots, f_{\ell-1})$, $K=(f_{1},
\ldots, f_{\ell+s})$, and consider the natural map of special fiber
rings $\varphi: \mathcal{F}(K) \longrightarrow \mathcal{F}(I)$. The
following are equivalent$:$

\begin{itemize}
\item[(i)] $I^{n}=(f_{1}, \ldots , f_{\ell-1})I^{n-1}+(f_{\ell},
\ldots, f_{\ell+s})^{n}$ for $n \gg 0$ ;

\item[(ii)]
${\rm{coker}}(\varphi)$ is a finite $\mathcal{F}(\f{a})$-module;

\item[(iii)] $\mathcal{F}(I)$ has embedding dimension at
most s locally at every minimal prime of dimension $\ell$.

\end{itemize}
\end{theorem}
\proof We apply Corollary \ref{finitecoker} and Proposition
\ref{maxidealpowers}.  \QED

\s
\s

\section{The Main Theorem}

\s

In this section we prove our main theorem about the core in
arbitrary characteristic. The proof uses reduction to the case of
equimultiple height one ideals, which we treat by means of the
results in the previous section. The reduction step on the other
hand requires the next two technical lemmas.

\begin{lemma}{\label{coker lifts}}
Let $R$ be a Noetherian local ring with infinite residue field $k$,
$I\/$ an $R$-ideal, and $J\/$ a reduction of $I\/$. Let $x$ be a
general element in $J\/$, write $x^{*}$ for the image of $x$ in
$[\mc{F}(I)]_{1}$, and let `${}^{\tratto}$' denote images in
$\overline{R}=R/(x)$.
\begin{itemize}

\item[(a)] The kernel of the natural map $\mc{F}(I)/x^{*}\mc{F}(I) \longrightarrow \mc{F}(\ol{I})$
is a finite-dimensional $k$-vector space.

\item [(b)] Let $\f{a}\subset K$ be $R$-ideals with $x \in \f{a}$ and $K \subset I$.
Consider the natural map of special fiber rings $\varphi:
\mathcal{F}(K) \longrightarrow \mathcal{F}(I)$ and write
$\ol{\varphi}\,$ for the induced map from $\mc{F}(\ol{K})$ to
$\mc{F}(\ol{I})$. Then ${\rm{coker}} \, (\ol{\varphi})$ is a finite
$\mathcal{F}(\ol{\f{a}})$-module if and only if ${\rm{coker}}
\,(\varphi)$ is a finite $\mathcal{F}(\f{a})$-module.
\end{itemize}
\end{lemma}
\proof To prove part (a) let $\mc{G}(I)$ and $\mc{G}(\ol{I})$ denote
the associated graded ring of $I\/$ and $\ol{I}$, respectively.
Consider the exact sequence
\[
0 \longrightarrow C \longrightarrow \mc{G}(I)/x^{*}\mc{G}(I)
\longrightarrow \mc{G}(\ol{I}) \longrightarrow 0.
\]
Since $x$ is general in $J$ and $J$ is a reduction of $I$, it
follows that $x$ is a superficial element of $I\/$. Thus $C$
vanishes in large degrees. Tensoring the above sequence with the
residue field $k$ we deduce that
\[
C\otimes_{R} k \longrightarrow \mc{F}(I)/x^{*}\mc{F}(I)
\longrightarrow \mc{F}(\ol{I}) \longrightarrow 0
\]
is exact and $C \otimes_{R} k$ is a finite-dimensional $k$-vector
space.

To prove part (b) notice that by (a) and the Snake Lemma the kernel
of the natural map
\[
{\rm{coker}} ( \varphi) /x^{*}{\rm{coker}} ( \varphi)
\twoheadrightarrow {\rm{coker}} (\ol{\varphi})
\]
is a finite-dimensional $k$-vector space as well. Hence
${\rm{coker}} ( \ol{\varphi})$ is finitely generated as a
$\mc{F}(\ol{\f{a}})$-module if and only if ${\rm{coker}}( \varphi)
/x^{*}{\rm{coker}} (\varphi)$ is finitely generated as a
$\mc{F}(\f{a})$-module. By the graded Nakayama Lemma the latter
condition means that ${\rm{coker}} (\varphi)$ is a finite
$\mc{F}(\f{a})$-module. \QED

\s
\s

The following two results use in an essential way the theory of
residual intersections. Let $R$ be a local Cohen-Macaulay ring,
$I\/$ an $R$-ideal, and $s$ an integer. Recall that ${\mathfrak a}
\, \colon I\/$ is a {\it geometric $s$-residual intersection} of
$I\/$ if ${\mathfrak a}$ is an $s$-generated $R$-ideal properly
contained in $I\/$ and ${\rm ht} \ {\mathfrak a} \, \colon I \geq s
$ as well as ${\rm ht} (I, {\mathfrak a} \, \colon I) \geq s+1$. The
ideal $I\/$ has the {\it Artin Nagata property} $AN_s^{-}$ if
$R/{\mathfrak a} \, \colon I\/$ is Cohen-Macaulay for every
geometric $i$-residual intersection ${\mathfrak a} \, \colon I$ and
every $i \leq s$.

\s

\begin{lemma}{\label{abars}}
Let $R$ be a local Cohen-Macaulay ring with infinite residue field
and assume that $R$ has a canonical module. Let $I\/$ be an
$R$-ideal with analytic spread $\ell >0$, and suppose that $I\/$
satisfies $G_{\ell}$ and ${\rm{AN}}_{\ell-1}^{-}$. Let $J\/$ be a
minimal reduction of $I\/$ and $K$ an $R$-ideal with $J \subset K
\subset I$. Consider the natural map of special fiber rings
$\varphi: \mathcal{F}(K) \longrightarrow \mathcal{F}(I)$. Assume
that ${\rm{coker}} \, \varphi$ has dimension at most $\ell-1$ as a
module over $\mc{F}(J)$. Write $\mc{A}=\mc{A}(J)$ for the set
consisting of all ideals $\f{a}$ such that $\f{a}:J$ is a geometric
$(\ell-1)$-residual intersection, $\mu(J/\f{a})=1$, and
${\rm{coker}}(\varphi)$ is a finite $\mathcal{F}(\f{a})$-module. For
$t$ a positive integer let $H$ be an $R$-ideal with ${\rm{ht}} (J,
J^{t}:H) \geq \ell$.

Then
\[
H \cap \underset{\f{a} \in \mathcal{A}} \bigcap ({J}^{t}, \f{a})
\subset {J}^{t}.
\]
\end{lemma}
\proof  We prove the lemma by induction on $\ell$. First let
$\ell=1$. As $I\/$ satisfies $G_1$, $J\/$ does and hence $0:J$ is a
geometric residual intersection of $J\/$. Thus $\mc{A}=\{0\}$ and
the assertion is clear. Therefore we may assume that $\ell \geq 2$.
Let $b \in H$ and suppose that $b \in J^{j-1} \backslash J^{j}$ for
some $j\/$ with $1 \leq j \leq t$. We are going to prove that there exists an
ideal $\f{a} \in \mc{A}$ with $b \notin (J^{t},\f{a})$. Since
$(J^{t}, \f{a})\subset J$ we may assume that $b \in J$.

We first reduce to the case where $I$ has positive height. Let
`${}^{\tratto}$' denote images in $\overline{R}=R/0:I$. Notice that
$0:I$ is a geometric 0-residual intersection of $I\/$ since $I\/$
satisfies $G_{1}$. Therefore $\ol{R}$ is Cohen-Macaulay by the
$AN_{0}^{-}$ condition, and ${\rm{ht}} \ \ol{I} >0$. Furthermore $I
\cap (0:I)=0$ according to \cite[1.7.c]{U94}. Thus the canonical
epimorphism $R \twoheadrightarrow \overline{R}\,$ induces
isomorphisms $I^{m} \simeq \ol{I}^{m}$  and $J^m \simeq
\overline{J}^m$ for every $m \geq 1$. Therefore $\overline{b} \in
\overline{J}^{j-1} \setminus \overline{J}^j$. Furthermore
$[\mc{G}(I)]_{m} \simeq [\mc{G}(\ol{I})]_{m}$  for $m \geq 1$, and
$\mc{F}(I) \simeq \mc{F}(\ol{I})$ for $m \geq 1$. Hence
$\ell(\ol{I})=\ell(I)$ and $\ol{J}$ is a minimal reduction of
$\ol{I}$. As ${\rm ht} \ 0 : I =0$ it follows that $\ol{I}$
satisfies $G_{\ell}$, and since $I \cap (0:I)=0$ the ideal $\ol{I}$
satisfies ${\rm{AN}}_{\ell-1}^{-}$ according to \cite[2.4.b]{JU}.
Obviously the cokernel of the induced map $\ol{\varphi}
:\mathcal{F}(\ol{K}) \longrightarrow \mathcal{F}(\ol{I})$ has
dimension at most $\ell-1$ as an $\mc{F}(\ol{J})$-module. Every
ideal in ${\mathcal A}(\overline{J}) $ is of the form
$\overline{{\mathfrak a}}$ for some ${\mathfrak a} \in {\mathcal A}
(J)$. Indeed if $\ol{\f{a}} \in \mc{A}(\ol{J})$ then there exists an
$\ell-1$-generated ideal $\f{a} \subset J$ whose image in $\ol{R}$
is $\ol{\f{a}}$. Since $J \cap (0:I)=0$ we have
$\f{a}:J=(0:I,\f{a}):J$, and it follows that $\f{a}:J$ is a
geometric $\ell-1$-residual intersection. Notice that a minimal
generating set of $\f{a}$ forms part of a minimal generating set of
$J\/$, hence $\mu(J/\f{a})=1$. Furthermore ${\rm{coker}} (\varphi)$
is a finite $\mc{F}(\f{a})$-module because $\mc{F}(I) \simeq
\mc{F}(\ol{I})$. Finally ${\rm ht} (\overline{J}, \overline{J}^t :
\overline{H}) \geq {\rm ht} (J, J^{t}:H) \geq \ell \,$ because ${\rm
ht} \ 0:I =0$. Therefore we may replace $R$ by $\ol{R}$ and assume
that ${\rm{ht}} \, I>0$. With this additional assumption we now
prove that $b \notin (J^{t}, \f{a})$ for some $\f{a} \in \mc{A}$.

Notice that ${\rm{ht}} \ J:I \geq \ell$ according to \cite[2.7]{JU}.
Since $I\/$ satisfies $G_{\ell}$ it then follows that $J\/$
satisfies $G_{\infty}$. Again as ${\rm{ht}} \ J:I \geq \ell$, the
property ${\rm{AN}}_{\ell-1}^{-}$ passes from $I\/$ to $J\/$ by
\cite[1.12]{U94}. Now $J\/$ satisfies the sliding depth condition
according to \cite[1.8.c]{U94}. In particular ${\rm Sym} (J/J^2)
\simeq \mc{G}(J)$ via the natural map and these algebras are
Cohen-Macaulay by \cite[6.1]{HSV}.

The proof of \cite[4.2]{PU} shows that $\ol{b} \notin \ol{J}^{j}$,
where now `${}^{\tratto}$' denotes images in $\overline{R}=R/(x)$
for a general element $x$ in $J\/$. By the general choice of $x$ in
$J$ and since $J$ is a reduction of $I$, we have $\ell(\ol{I}) \leq
\ell(I)-1$ and then Lemma~\ref{coker lifts}(a) shows that
$\ell(\ol{I})=\ell(I)-1$. Again because $x$ is a general element and
${\rm ht} \, J > 0$, it follows that $x$ is $R$-regular. For the same
reasons and because $\mc{G}(J)$ is Cohen-Macaulay, the leading form
$x^*$ of $x$ in $\mc{G}(J)$ is regular on $\mc{G}(J)$, which gives
 $\mc{G}(J)/x^*\mc{G}(J)\simeq \mc{G}(\ol{J})$.
Therefore ${\rm Sym} (\ol{J}/\ol{J}^2) \simeq \mc{G}(\ol{J})$,
forcing $\ol{J}$ to satisfy $G_{\infty}$. Hence $\ol{I}$ satisfies
$G_{\ell-1}$, because ${\rm{ht}} \ \ol{J}:\ol{I} \geq \ell -1$. As
$x$ is an $R$-regular element it is easy to see that $\ol{I}$ is
${\rm{AN}}_{\ell-2}^{-}$.

Again by the general choice of $x$ the cokernel of the natural map
from $\mc{F}(\ol{K})$ to $\mc{F}(\ol{I})$ has dimension at most
$\ell-2$ as a $\mc{F}(\ol{J})$-module. Finally ${\rm ht}
(\overline{J}, \overline{J}^t : \overline{H}) \geq \ell -1$ and
according to Lemma~\ref{coker lifts}(b) every ideal of
$\mc{A}(\ol{J})$ is of the form $\ol{\f{a}}$ for some $\f{a} \in
\mc{A}$. Thus by the induction hypothesis, $\ol{b} \notin
(\ol{J}^{t}, \ol{\f{a}})$ for some $\f{a} \in \mc{A}$. Hence $b
\notin (J^{t}, \f{a})$. \QED

\s

We are now ready to prove our main result.

\begin{theorem} \label{emdimformula}
Let $R$ be a local Cohen-Macaulay ring with infinite perfect residue
field, let $I\/$ be an $R$-ideal with $g={\rm{ht}} \, I >0$ and
$\ell=\ell(I)$, and let $J\/$ be a minimal reduction of $I\/$ with
$r=r_{J}(I)$. Suppose that the special fiber ring $\mathcal{F}(I)$
of $I\/$ has embedding dimension at most one locally at every
minimal prime of dimension $\ell$. If $\ell > g$ further assume that
$R$ is Gorenstein, $I\/$ satisfies $G_{\ell}$ and ${\rm{depth}} \
R/I^{j} \geq {\rm{dim}} \ R/I-j+1$ for $1\leq j \leq \ell-g$. Then
\[
{\rm{core}}(I)= J^{n+1}:I^{n}
\]
for every $n \geq {\rm{max}} \{r-\ell+g,0 \}$.

\end{theorem}

\proof  According to \cite[2.9(a)]{U94} the ideals $I\/$ and
$I\hat{R}$ satisfy ${\rm{AN}}_{\ell-1}^{-}$, and hence are
universally weakly $\ell -1$ residually $S_2$ in the sense of
\cite[p. 203]{CEU}. Therefore \cite[4.8]{CPU01} shows that
$\core{I}\hat{R}=\core{I\hat{R}}$. Thus we may pass to the
completion of $R$ and assume that $R$ has a canonical module. Let
$f_{1}, \ldots , f_{\ell+1}$ be general elements in $I\/$. The ideal
$J^{n+1}:I^{n}$ for $n\geq {\rm{max}} \{r-\ell+g,0 \}$ is
independent of the minimal reduction $J\/$ and of $n$, as can be
seen from \cite[5.1.6]{HySm1} if $\ell =g$ and from \cite[2.3]{PU}
if $\ell > g$. Hence we may assume that $J=(f_{1}, \ldots ,
f_{\ell})$ and $n \gg 0$. We use the notation of Lemma \ref{abars}
with $K=(J, f_{\ell+1})$, $t=n+1$, and $H$ the intersection of all
primary components of $J^{n+1}$ of height $< \, \ell$. Notice that
${\rm{coker}} (\varphi)$ has dimension at most $\ell-1$ as a module
over $\mc{F}(J)$ according to Theorem~\ref{I^n decomposition}. Hence
the assumptions of Lemma~\ref{abars} are satisfied.

Let $\f{a} \in \mc{A}$ be as in Lemma~\ref{abars}. Write
`${}^{\tratto}$' for images in $\overline{R}=R/\f{a}:I$. Notice that
$\overline{R}$ is Cohen-Macaulay and by \cite[1.7(a)]{U94},
${\rm{ht}} \ \ol{I}>0$. Hence ${\rm{ht}} \
{\overline{I}}=\ell(\overline{I})=1$. Now \cite[3.4]{PU} shows that
${\rm{core}} (\ol{I})= \ol{J}^{n+1}: \underset{y \in \ol{I}}
\sum(\ol{J},y)^n$. Notice that $\ol{K}^{n} \subset \underset{y \in
\ol{I}} \sum(\ol{J},y)^n \subset \ol{I}^{n}$ and that
$\ol{K}^{n}=\ol{I}^{n}$ according to Corollary ~\ref{finitecoker}.
Hence
${\rm{core}}(\overline{I})={\overline{J}}^{n+1}:{\overline{I}}^{n}$.

On the other hand by \cite[4.5]{CPU01},
${\rm{core}}(\overline{I})=(\overline{\alpha_{1}}) \cap \ldots \cap
(\overline{\alpha_{\gamma}})$ for some integer $\gamma$ and $\gamma$
general principal ideals $(\overline{\alpha_{1}}), \ldots ,
(\overline{\alpha_{\gamma}})$ in $\ol{I}\/$. Notice that
$(\f{a},\alpha_{i})$ are reductions of $I\/$, hence ${\rm{core}}(I)
\subset \overset{\gamma}{\underset{i=1} \bigcap}
(\f{a},\alpha_{i})$. Therefore $\overline{{\rm{core}}(I)} \subset
\overline{\overset{\gamma}{\underset{i=1} \bigcap}
(\f{a},\alpha_{i})} \subset
\overset{\gamma}{\underset{i=1}{\bigcap}} (\overline{\alpha_{i}})=
{\rm{core}}(\overline{I})$. As ${\rm{core}}(\overline{I}) =
{\overline{J}}^{n+1}:{\overline{I}}^{n}$ we obtain
\begin{align*}
            {\rm{core}}(I) &\subset (J^{n+1}, \f{a}:I):I^{n}\\
                    &=(J^{n+1}, (\f{a}:I)\cap I):I^{n}\\
                      &=(J^{n+1},\f{a}):I^{n}.
\end{align*}
The last equality holds because $(\f{a}:I)\cap I=\f{a}$ by
\cite[1.7(c)]{U94}. It follows that
\begin{equation}\label{(1.4)}
{\rm{core}}(I) \subset  \underset{\f{a} \in \mathcal{A}}{\bigcap}
(J^{n+1},\f{a}):I^{n}.
\end{equation}

Next we show that
\begin{equation}\label{(1.5)}
 {\rm{core}}(I) \subset H:I^{n},
\end{equation}
or equivalently $({\rm{core}}(I))_{\f{p}} \subset (H:I^{n})_{\f{p}}$
for every prime ideal $\f{p}$ with ${\rm{dim}}\, R_{\f{p}} < \ell$.
Indeed by \cite[2.7]{JU}, $J_{\f{p}}=I_{\f{p}}$, and hence
$J_{\f{p}}^{n}=I_{\f{p}}^{n}$. Thus $({\rm{core}}(I))_{\f{p}}
\subset J_{\f{p}} \subset J_{\f{p}}^{n+1}:J_{\f{p}}^{n}= H_{\f{p}}:
I_{\f{p}}^{n}$.

Finally $J^{n+1}:I^{n} \subset {\rm{core}}(I)$, as can be seen from
the proof of \cite[4.5]{PU} via \cite[5.1.6]{HySm1} if $\ell=g$ and
from \cite[4.8]{PU} otherwise. Hence (\ref{(1.4)}), (\ref{(1.5)})
and Lemma \ref{abars} imply that
\begin{align*}
        J^{n+1}:I^{n} \subset  {\rm{core}}(I) & \subset (H \cap \underset{\f{a} \in
          \mathcal{A}} \bigcap (J^{n+1}, \f{a})):I^{n}\\
          &\subset J^{n+1}:I^{n}
\end{align*}
Therefore ${\rm{core}}(I)= J^{n+1}:I^{n}$. \QED

\s
\s

\section{Applications}

\s

In this section, we collect several instances where the assumption
on the generic embedding dimension of the special fiber ring
required in Theorem~\ref{emdimformula} holds automatically.

\s

\begin{theorem}\label{cor integral}
Let $k$ be an infinite field, $R'$ a positively graded geometrically
reduced Cohen-Macaulay $k$-algebra, and $R$ the localization of $R'$
at the homogeneous maximal ideal. Let $K$ be an $R$-ideal generated
by forms in $R'$ of the same degree, let $f$ be an element of $R$
integral over $K$, and write $I=(K,f)$. Set $g={\rm ht} \, {I}
>0$, $\ell=\ell(I)$ and let $J\/$ be a minimal reduction of $\, I \/$
with $r=r_{J}(I)$. If $\ell > g$ suppose that $R'$ is Gorenstein,
$I\/$ satisfies $G_{\ell}$ and ${\rm{depth}} \ R/I^{j} \geq
{\rm{dim}} \ R/I-j+1$ for $1\leq j \leq \ell-g$. Then
\[
{\rm{core}}(I)= J^{n+1}:I^{n}
\]
for every $n \geq {\rm{max}} \{r-\ell+g,0 \}$.
\end{theorem}

\proof
Observe that
\[
R=R'_{R'_{+}} \hookrightarrow
S=(R'\otimes_{k}\ol{k})_{(R'\otimes_{k}\ol{k})_{+}}
\]
is a flat local extension. Furthermore according to
\cite[2.9(a)]{U94} the ideals $I\/$ and $IS$ are universally weakly
$\ell -1$ residually $S_2$.
Therefore
\cite[4.8]{CPU01} shows that $\core{I}S=\core{IS}$. Thus, replacing
$k$ by $\ol{k}$ we may suppose that $k$ is perfect and $R'$ is
reduced.

Write $K=(f_1, \ldots, f_{m})$ where $f_{1}, \ldots, f_{m}$ are
forms of the same degree. Now $\mathcal{F}(K)\simeq k[f_{1}, \ldots
, f_{m}]$ is a subalgebra of $R'$ and thus is reduced. Let $\f{p}$
be a minimal prime of $\mathcal{F}(I)$ of dimension $\ell$ and write
$\f{q}$ for its contraction to $\mathcal{F}(K)$. As $K$ is a
reduction of $I\/$ we have $\ell(K)=\ell(I)$ and therefore
${\rm{dim}} \, \mathcal{F}(K)=\ell={\rm{dim}} \, \mathcal{F}(I)$.
Furthermore $\mathcal{F}(I)$ is finitely generated as a module over
$\mathcal{F}(K)$. It follows that $\f{q}$ is a minimal prime of
$\mathcal{F}(K)$. Since $\mathcal{F}(K)$ is reduced the localization
${\mathcal{F}(K)}_{q}$ is a field, say $L$. Now
$\mathcal{F}(I)_{\f{p}}$ is a localization of an $L$-algebra
generated by a single element, namely the image of $f$. Hence
$\mathcal{F}(I)_{\f{p}}$ has embedding dimension at most one. As
this holds for every minimal prime $\f{p}$ of dimension $\ell$, the
result follows from Theorem~\ref{emdimformula}. \QED

\s

\begin{remark} Notice that taking $f=0$ in Theorem~\ref{cor integral} we obtain
Theorem~\ref{formsofthesamedegree} of the Introduction. There is a
graded and a global version of the latter theorem if the ideal $I$
is zero-dimensional. Thus, let $I'$ be a homogeneous $R'$-ideal with
$I'R=I$ and let $J'$ be an $R'$-ideal generated by ${\rm{dim}} \,  R'$ general
$k$-linear combinations of homogeneous minimal generators of $I'$.
One has
\[
{\rm gradedcore}(I')=\core{I'}= J'^{\, n+1}:I'^{\, n}  \hspace{1cm}
\mbox{ for every \ } n \geq r
\]
where ${\rm gradedcore}(I')$ stands for the intersection of all
homogeneous reductions of $I'$.

In fact, since $I'$ is zero-dimensional and generated by forms of
the same degree, the first equality obtains by \cite[4.5]{CPU01} and
\cite[2.1]{PUV}, whereas the second equality follows from
Theorem~\ref{formsofthesamedegree} and \cite[2.1]{PUV}.
\end{remark}

\s

\begin{theorem}\label{regular}Let $R$ be a local
Cohen-Macaulay ring with infinite perfect residue field, let $I\/$
be an $R$-ideal with $g={\rm{ht}} \, {I} >0$ and $\ell=\ell(I)$, and let
$J\/$ be a minimal reduction of $ \, I\/$ with $r=r_{J}(I)$. Suppose
that $I=(K,f)$, where $\ell(I) \geq \ell(K)$ and the special fiber
ring $\mathcal{F}(K)$ satisfies Serre's condition $R_{1}$. If $\ell
>g$ further assume that $R$ is Gorenstein, $I\/$ satisfies $G_{\ell}$,
and ${\rm{depth}} \ R/I^{j} \geq {\rm{dim}} \ R/I-j+1$ for $1\leq
j \leq \ell-g$. Then
\[
{\rm{core}}(I)= J^{n+1}:I^{n}
\]
for every $n \geq {\rm{max}} \{r-\ell+g,0 \}$.

\end{theorem}

\proof According to Theorem~\ref{emdimformula} it suffices to prove
that $\mc{F} (I)$ has embedding dimension at most one locally at
every minimal prime of dimension $\ell = {\rm{dim}} \, \mc{F}(I)$.

Let $A=\mathcal{F}(K)$ and $B=\mathcal{F}(I)$. Let $\f{q}$ be a
prime ideal of $B$ of dimension $\ell$ and write $\f{p}=\f{q} \cap
A$. We claim that ${\rm{dim}} \, A_{\f{p}} \leq 1$. The affine
domain $B/\f{q}$ is generated by one element as an algebra over
$A/\f{p}$. Therefore
\[
{\rm{dim}} \ A/\f{p} \geq {\rm{dim}} \ B/\f{q} -1= {\rm{dim}} \, B -1.
\]
Hence
\[
{\rm{dim}} \, A_{\f{p}} \leq {\rm{dim}} \, A -{\rm{dim}} \ A/\f{p}
\leq {\rm{dim}} \, A - {\rm{dim}} \, B +1 =\ell(K) -\ell +1 \leq 1.
\]

Since ${\rm{dim}} \, A_{\f{p}} \leq 1$ our assumption gives that
$A_{\f{p}}$ is regular. Now we consider the exact sequence of
modules of differentials,
\[
B_{\f{q}} \otimes_{A_{\f{p}}} \Omega_{k}(A_{\f{p}})
 \longrightarrow \Omega_{k}(B_{\f{q}})
 \longrightarrow  \Omega_{A_{\f{p}}}(B_{\f{q}}) \longrightarrow 0.
\]

As $A_{\f{p}}$ is regular  and $k$ is perfect it follows that
$\mu_{A_{\f{p}}}(\Omega_{k}(A_{\f{p}})) \leq {\rm{dim}} \,
A_{\f{p}}+ {\rm{trdeg}}_{k} \, A/\f{p}$. Hence $\mu_{B_{\f{q}}}(
B_{\f{q}} \otimes_{A_{\f{p}}} \Omega_{k}(A_{\f{p}})) \leq {\rm{dim}}
\, A_{\f{p}}+ {\rm{trdeg}}_{k} \, A/\f{p}$. Since $B$ is generated by
one element as an $A$-algebra, the $B_{\f{q}}$-module
$\Omega_{A_{\f{p}}}(B_{\f{q}})$ is cyclic. Computing numbers of
generators along the above exact sequence we obtain
\[
\mu_{_{B_{\f{q}}}}(\Omega_{k}(B_{\f{q}})) \leq {\rm{dim}} \,
A_{\f{p}}+{\rm{trdeg}}_{k} \, A/\f{p}+1.
\]
On the other hand by \cite[Satz 1$($a$)$]{BK},
\[
\mu_{_{B_{\f{q}}}}(\Omega_{k}(B_{\f{q}}))={\rm{edim}} \, B_{\f{q}}+{\rm{trdeg}}_{k} \, B/\f{q}.
\]
We conclude that
\begin{align*} {\rm{edim}} \, B_{\f{q}} &\leq {\rm{dim}} \,
A_{\f{p}}+{\rm{trdeg}}_{k} \, A/\f{p}-{\rm{trdeg}}_{k} \, B/\f{q}+1 \\
&= {\rm{dim}} \, A_{\f{p}}+{\rm{dim}} \ A/{\f{p}}-{\rm{dim}} \, B+1\\
 &\leq {\rm{dim}} \, A-{\rm{dim}} \, B +1\\
 &\leq 1.
\end{align*}
\QED

\s

The above result can be considered as a generalization of the case
of second analytic deviation one treated in \cite[4.8]{PU}. In this
case the minimal number of generators of $I\/$ exceeds $\ell$ by at
most one, and we can choose $K$ to be $J\/$ in
Theorem~\ref{regular}. But then $\ell(K)=\ell$ and $\mc{F} (K)$
satisfies $R_1$, being a polynomial ring over $k$.

Also observe that the condition $\ell \geq \ell(K)$ in
Theorem~\ref{regular} is always satisfied if $I\/$ is primary to the
maximal ideal. Here is another situation where this inequality holds
automatically:

\s

\begin{remark}\label{gradedregular}
Let $k$ be an infinite perfect field, $R'$ a positively graded
Cohen-Macaulay $k$-algebra, and $R$ the localization of $R'$ at the
homogeneous maximal ideal. Let $K$ be an $R$-ideal generated by
forms in $R'$ of the same degree $e$, let $f$ be a form in $R'$ of
degree at least $e$, write $I=(K,f)$, and assume that the subalgebra
$k[K_e]$ of $R'$ satisfies Serre's condition $R_{1}$. Set $g={\rm{ht}}
\, I >0$, $\ell=\ell(I)$ and let $J\/$ be a minimal reduction of $I\/$
with $r=r_{J}(I)$. If $\ell > g$ further suppose that $R'$ is
Gorenstein, $I\/$ satisfies $G_{\ell}$ and ${\rm{depth}} \ R/I^{j}
\geq {\rm{dim}} \ R/I-j+1$ for $1\leq j \leq \ell-g$. Then
\[
{\rm{core}}(I)= J^{n+1}:I^{n}
\]
for every $n \geq {\rm{max}} \{r-\ell+g,0 \}$.

\end{remark}

\proof After rescaling the grading we can identify the subalgebra
$k[K_e]$ of $R'$ with $\mathcal{F}(K)$, which shows that the latter
ring satisfies $R_1$. Thus to apply Theorem~\ref{regular} it
suffices to verify that $\ell(I) \geq \ell(K)$.  Comparing Hilbert
functions it follows that $\ell(I)={\rm dim} \, \mathcal{F}(I) \geq
{\rm{dim}} \,\mathcal{F}(K) = \ell(K)$ once we have proved the
injectivity of the natural map $\varphi: \mathcal{F}(K)
\longrightarrow \mathcal{F}(I)$. To show the latter, write $\f{m}$
for the maximal ideal of $R$. Let $F$ be a form of degree $s$ in
$\mc{F}(K)$ such that $\varphi(F)=0$. Then $F \in \f{m}I^{s}$ as an
element of $R'\subset R$. In $R'$ the form $F$ has degree $se$,
whereas the nonzero homogeneous elements of $\f{m}I^{s}$ have
degrees at least $se+1$. Therefore $F=0$.

\QED


\s

\section{Examples}

\ms

In this section we present several examples showing that the various
assumptions in our theorems are in fact necessary. We will always
use zero-dimensional ideals in local Gorenstein rings, so that the
property $G_{\ell}$ as well as the depth conditions for the powers
of the ideal hold automatically.

\s

The first example illustrates that
Theorem~\ref{formsofthesamedegree} is no longer true if the ring
$R'$ fails to be geometrically reduced, even if it is a domain and
all the other assumptions of the theorem are satisfied.

\begin{example}\label{reduced}
{\rm Let $k_0$ be a field of characteristic $p > 0$ and let
$k=k_0(s,t)$ be the rational function field in two variables.
Consider the ring $R'=k[x,y,z]/(x^p-sz^p, y^p-tz^p)$.  This ring is
a one-dimensional standard graded Gorenstein domain. Indeed, the
elements $x^p, y^p, z^p$ generate an ideal of grade $3$ in
$k_0[x,y,z]$, and $x^p-sz^p, y^p-tz^p$ are obtained from generic
linear combinations of these elements by localization, change of
variables, and descent. Thus \cite[{\rm Theorem~(b)}]{Ho} shows that
$x^p-sz^p, y^p-tz^p$ generate a prime ideal in the ring $k[x,y,z]$.

On the other hand, $R'$ is not geometrically reduced. Indeed, after
tensoring with the algebraic closure $\ol{k}$ of $k$ we obtain
\[R' \otimes_{k} \ol{k} \simeq
\ol{k}[x,y,z]/((x-z\sqrt[p]{s})^{p},(y-z\sqrt[p]{t})^{p}),\] which
is not a reduced ring.

Let $(R, \f{m})$ denote the localization of $R'$ at the homogenous
maximal ideal. We claim that
\[\core{\f{m}} \not= J^{n+1} \colon
\f{m}^{n} \hspace{.15in} \mbox{for every $n \gg 0$ and every minimal
reduction $J\/$ of $\f{m}$}.\] Indeed, the $p^{\rm th}$ power of any
general linear form in $R$ generates the ideal $Rz^p$. As the core
of $\f{m}$ is a finite intersection of principal ideals generated by
general linear forms \cite[4.5]{CPU01}, it follows that $Rz^p
\subset \core{\f{m}}$. On the other hand, according to
\cite[3.2$($a$)$]{PU} one has $J^{n+1} \colon \f{m}^{n}= Rz^{n+1}
\colon \f{m}^{n}$, since $Rz$ is a minimal reduction of $\f{m}$. As
$R'/R'z^{n+1}$ is a standard graded Artinian Gorenstein ring with
$a$-invariant $n+2p-2$ it follows that $Rz^{n+1} \colon \f{m}^{n}=
\f{m}^{2p-1}$, which does not contain $Rz^p$. }
\end{example}

\s

The next example shows that the assumption in
Theorem~\ref{emdimformula} on the local embedding dimension of the
special fiber ring is sharp: If we allow the local embedding
dimension to be $2$ the statement of the theorem is no longer true
even in the presence of the other conditions.

\begin{example}\label{emb2}
{\rm Let $k$ be an infinite perfect field of characteristic $2$,
$R=k[x,y]_{(x,y)}$ a localized polynomial ring and $I=(x^6, x^5y^3,
x^4y^4, x^2y^8, y^9)$. Using Macaulay 2 \cite{M2} one computes the
special fiber ring of $I\/$ to be
\[ \mc{F}(I) \simeq k[a,b,c,d,e]/(b^2, bd, cd, d^2, c^2 - ad),\]
where $a,b,c,d,e$ are variables over $k$.

This ring has a unique minimal prime ideal $\p$, which is generated
by the images of $b, c, d$, and one easily sees that $[\mc{F}
(I)]_{\p}$ has embedding dimension $2$.

We claim that
\[\core{I} \not= J^{n+1} \colon
I^{n} \hspace{.15in} \mbox{for every $n \gg 0$ and every minimal
reduction $J\/$ of $I\/$}.\]

Indeed, $H=(x^6,y^9)$ is a minimal reduction of $I\/$ with
$r_H(I)=2$. Thus $J^{n+1}: I^n=H^3:I^2$ according to \cite[2.3]{PU}.
On the other hand, using the algorithm of \cite[3.6]{PUV} it has
been shown in \cite[3.9]{PUV} that $\core{I} \not= H^3:I^2$.}
\end{example}

\s

The next example shows that in Theorem~\ref{cor integral} and
Theorem~\ref{regular} it is essential to assume that either $f$ is
integral over $K$ or else the special fiber ring $\mc{F} (K)$
satisfies Serre's condition $R_{1}$.

\begin{example}\label{serre}
{\rm Let $k$ be an infinite perfect field of characteristic $2$,
$R=k[x,y]_{(x,y)}$ a localized polynomial ring, $K =(x^9, x^5y^4,
x^3y^6, x^2y^7)$, which is an ideal generated by monomials of the
same degree, $f=y^8$ and $I=(K, f)$.

Again we claim that
\[\core{I} \not= J^{n+1} \colon
I^{n} \hspace{.15in} \mbox{for every $n \gg 0$ and every minimal
reduction $J\/$ of $I\/$}.\]

The ideal $H=(x^9,y^8)$ is a minimal reduction of $I\/$ with
$r_H(I)=2$, and hence $J^{n+1}: I^n=H^3:I^2$ by \cite[2.3]{PU}. On
the other hand, using the algorithm of \cite[3.6]{PUV} and Macaulay
2 {\cite{M2}} we can compute \[\core{I}= H^3:I^2 + (xy^{12},y^{13})
\varsupsetneq H^3 :I^2.\]}
\end{example}

\s
\s
\s

\s

\end{document}